\documentclass[10pt]{article}

 \usepackage{authblk}
\newtheorem{rem}{Remark}[section]
\usepackage{amsmath,amssymb} % estetics
\usepackage{graphicx}        % standard LaTeX graphics tool
\usepackage{multicol}        % used for the two-column index
\usepackage[bottom]{footmisc}% places footnotes at page bottom

\newcommand{\lga}[2]{\langle #1,#2\rangle_\Gamma}
\newcommand{\bbR}{\mathbb{R}}
\newcommand{\mcK}{\mathcal{K}}
\newcommand{\bfn}{\boldsymbol{n}}
\newcommand{\bfM}{\boldsymbol{M}}
\newcommand{\bft}{\boldsymbol{t}}
\newcommand{\bfx}{\boldsymbol{x}}
\newcommand{\bfkappa}{\boldsymbol{\kappa}}
\newcommand{\bfI}{\boldsymbol{I}}
\newcommand{\strain}{\boldsymbol{\varepsilon}}
\newcommand{\stress}{\boldsymbol{\sigma}}
\newcommand{\mcCP}{D}
\newcommand{\Div}{\text{\bf div}}       % 

\begin{document}

\title{\bf Augmented Lagrangian Method for Thin Plates with Signorini Boundaries}
\date{}

% defining authors
\author[$\star$]{Erik Burman}
\author[$\dagger$]{Peter Hansbo}
\author[$\ddagger$]{Mats G. Larson}
\affil[$\star$]{\footnotesize  Department of Mathematics, University College London, London, UK--WC1E  6BT, United Kingdom}

\affil[$\dagger$]{\footnotesize Department of Mechanical Engineering, J\"onk\"oping University, SE-55111 J\"onk\"oping, Sweden}

\affil[$\ddagger$]{\footnotesize Department of Mathematics and Mathematical Statistics, Ume{\aa} University, SE-90187 Ume{\aa}, Sweden}

\maketitle

\begin{abstract}
We consider $C^1$-continuous approximations of the Kirchhoff plate problem in combination
with a mesh dependent augmented Lagrangian method on a simply supported Signorini boundary.
\end{abstract}

\section{Introduction}\label{sec:1}
%\subsection{The Augmented Lagrangian Method for a model Signorini problem\label{sec:1}}

To introduce the augmented Lagrangian method we first consider a simple Poisson problem, find $u$ such that
\begin{equation}\label{poiss}
-\Delta u= f ~\text{in} ~\Omega ,\qquad u=g ~ \text{on} ~\Gamma
\end{equation}
where $\Omega$ is a bounded domain with boundary $\Gamma := \partial \Omega$ 
and exterior unit normal $\boldsymbol{n}$,

The Lagrange multiplier approach to prescribing $u=g$ is to seek stationary points to
\begin{equation}
\mathcal{L}(v,\mu) := \frac12 a(v,v)  -
\langle{\mu},{v-g}\rangle_\Gamma  - ({f},{v}) 
\end{equation}
where 
\begin{equation}
({f},{v}) := \int_{\Omega}f v \, d\Omega, \; a(u,v) := \int_\Omega \nabla u\cdot\nabla v\, d\Omega
\end{equation}
and $\langle{\cdot},{\cdot}\rangle_\Gamma$ denotes the
$H^{-1/2}/H^{1/2}$-duality pairing. Whenever the arguments are smooth
enough we define,
\begin{equation}
 \langle{\mu},{v-g}\rangle_\Gamma := \int_{\Gamma}\mu (v-g) \, ds 
\end{equation}
{Stationary} points are given by finding $(u,\lambda)\in H^1(\Omega)\times H^{-1/2}(\Gamma)$ such that
\begin{equation}
a(u,v)-\langle{\lambda},{v}\rangle_\Gamma = (f,v)\quad \forall v\in H^1(\Omega)
\end{equation}
\begin{equation}\label{eq:weakform2}
\langle{\mu},{u}\rangle_\Gamma = \langle{\mu},{g}\rangle_\Gamma\quad \forall \mu\in H^{-1/2}(\Gamma) 
\end{equation}
{Formally,} the Lagrange multiplier is given by $\lambda = \partial_n u$, where $\partial_n{v}:={\boldsymbol n}\cdot\nabla v$. In a discretization of this problem, the approximation of the multiplier and the displacement must fulfil an {\em inf--sup}\/ condition ensuring that the problem will not be overconstrained.

We now augment the Lagrangian by a penalty term and seek stationary points to
\begin{equation}\label{eq:aug1}
\mathcal{L}(v,\mu) := \frac12 a(v,v)  -
\lga{\mu}{v-g}  +\frac12\|\gamma^{1/2} (v-g)\|^2_\Gamma - ({f},{v}) 
\end{equation}
leading to the problem of finding $(u,\lambda)\in H^1(\Omega)\times H^{-1/2}(\Gamma)$ such that
\begin{equation}
 a(u,v)  -
\lga{\lambda}{v}  +\lga{\gamma\, u}{v} -\lga{\mu}{u} =  ({f},{v})+\lga{\gamma\, g}{v}-\lga{\mu}{g}  \\
\end{equation}
for all $(v,\mu)\in H^1(\Omega)\times H^{-1/2}(\Gamma)$. The
discretization of this problem requires the same careful balance
between approximation spaces for the primal variable and the
multiplier as does the standard Lagrange multiplier method. Indeed if
we introduce the space
\begin{equation}
V_h := \{v_h \in H^1(\Omega): v_h\vert_K \in \mathbb{P}_k(K), \,
\forall K \in \mathcal{T}_h \},\quad \mbox{ for } k \ge 1
\end{equation}
where $\mathcal{T}_h$ is a conforming quasi-uniform partition of
$\Omega$ and $\mathbb{P}_k(K)$ denotes the set of polynomials of
degree less than or equal to $k$ on the element $K$ for the
discretization of $u$, we must find a mulitplier space $\Lambda_h$ such 
that the inf-sup condition is satisfied.
However, If we seek $u_h \in V_h$ with the discrete multiplier $\lambda_h := \partial_n u_h$, we recover Nitsche's method: 
\begin{equation}
 a(u_h,v)  -
\lga{\partial_n u_h}{v} -\lga{\partial_n v}{u_h} +\lga{\gamma\, u_h}{v}  =  ({f},{v})+\lga{g}{\gamma v-\partial_n v}  \\
\end{equation}
for all $v\in V_h$, with $\mu = \partial_n v$, which is stable with
the choice $\gamma = \gamma_0/h$, where $h$ is the local meshsize and
$\gamma_0$ large enough.

If we alternatively consider a stable discretization $\lambda_h\in \Lambda_h$, the discrete problem can be seen as seeking stationary points 
$(u_h,\lambda_h)\in V_h\times \Lambda_h$ to the modified Lagrangian
\begin{equation}
\mathcal{L}_h(v,\mu) := \frac12 a(v,v)    +\frac12\|\gamma^{1/2} (v-g-\gamma^{-1}\mu)\|^2_\Gamma -
\|\gamma^{-1/2}\mu\|_{\Gamma}^2- (f,v)_{\Omega} 
\end{equation}
which is obtained from $\mathcal{L}$ in (\ref{eq:aug1}) by rearranging terms and noting that the discrete multiplier is in $L_2(\Gamma)$.

We now turn to an inequality constraint on the boundary: $u \leq g$ on $\Gamma$. The corresponding Kuhn--Tucker conditions read:
\begin{equation}\label{eq:Kuhn1}
u -g \leq 0,\quad \lambda \leq 0, \quad \lambda (u-g)=0 \quad \text{on $\Gamma$}.
\end{equation}
These conditions can alternatively be written (cf. \cite{ChHi12})
\begin{equation}\label{eq:lambda1}
\lambda = -{\gamma}\,[u-g-\gamma^{-1}\, \lambda]_+ 
\end{equation}
where $\gamma \in \mathbb{R}^+$, $[x]_{+}=\max (x,0)$.
We can now, following  Alart and Curnier \cite{AlCu91}, define the following discrete augmented Lagrangian
\begin{equation}
\mathcal{L}_h(v,\mu) := \frac12 a(v,v)  +
\frac12 \Vert \gamma^{1/2}[v-g-\gamma^{-1} \mu]_+\Vert_\Gamma^2 -\|\gamma^{-1/2}\mu\|_{\Gamma}^2- ({f},{v})
\end{equation}
The corresponding Euler-Lagrange equations read: find $(u_h,\lambda_h)\in V_h\times \Lambda_h$ such that
\begin{equation}\label{aug1hb}
a({u_h},{v}) +\lga{{\gamma}\left[{u_h}-g-\gamma^{-1} \lambda_h\right]_+}{v} = ({f},{v}) \quad \forall {v}\in V_h
\end{equation}
and
\begin{equation}\label{aug2hb}
\lga{{\gamma}\left[{u_h}-g-\gamma^{-1}\lambda_h\right]_++  \lambda_h}{ \gamma^{-1}\mu} =0  \quad \forall \mu\in \Lambda_h 
\end{equation}
If $[{u_h}-g-\gamma^{-1}\lambda_h]_+=0$ (no contact) then $\lambda_h
=0$ and if $[{u_h}-g-\gamma^{-1}\lambda_h]_+ > 0$ (contact) we recover
the standard augmented formulation for the imposition of the Dirichlet
condition $u=g$. The multiplier approach \eqref{aug1hb}-\eqref{aug2hb}
using a stable pair $V_h \times \Lambda_h$
was shown to produce approximations of optimal accuracy in \cite{BuHaLa19}.
Now set $\lambda_h = \partial_n u_h$, $\mu = \partial_n v$ and seek $u_h\in V_h$ such that
\begin{equation}\label{eq:contact_Nitsche}
a(u_h,v)+\lga{{\gamma}\,[u_h-g-\gamma^{-1}\, \partial_n u_h]_+}{v-\gamma^{-1}\, \partial_n v} - \lga{\gamma^{-1}\,\partial_n u_h}{\partial_n v}= (f,v)
\end{equation}
for all $v\in  V_h$. With the choice $\gamma = \gamma_0/h$ this is the Nitsche method for Signorini problems first proposed in the context of elastic contact by
Chouly and Hild \cite{ChHi12}. For more information on augmented Lagrangian methods and variants thereof, see \cite{BuHa17}.

\section{The Kirchhoff plate model}
\label{subsec:2}
We now proceed formally to extend the discussion to the Kirchhoff plate model, posed on a domain 
$\Omega\subset \bbR^2$ with boundary $\Gamma = \partial\Omega$ and exterior unit normal 
$\bfn$. We seek an out--of--plane 
(scalar) displacement $u$ to which we associate the strain (curvature) tensor
\begin{equation}
\bfkappa (u) := \strain(\nabla u) := \frac12\left(\nabla\otimes (\nabla u) + (\nabla u ) \otimes \nabla \right) 
 = \nabla \otimes \nabla u%  = \nabla^2 u   
\end{equation}
and the plate stress (moment) tensor
\begin{align}\label{eq:plate-stress-tensor}
\bfM(u) :=\stress (\nabla u) := &\mcCP \left(\strain(\nabla u) + \nu (1- {\nu })^{-1}
\text{div}\nabla u \, \bfI \right)
\\
= &\mcCP \left( \bfkappa( u)+ \nu (1-\nu)^{-1} \Delta u \bfI \right)
\end{align}
where 
\begin{equation}\label{eq:mcCP}
\mcCP =  \frac{E t^3}{12(1+\nu)} 
\end{equation}
with $E$ the Young's modulus, $\nu$ the Poisson's ratio, and $t$ the 
plate thickness. We will use the standard convention that all quantities are positive downwards.

The Kirch\-hoff equilibrium problem takes the form: given the out--of--plane 
load (per unit area) $f$, find the displacement $u$ such that
\begin{align}
\text{div}\, \Div\, \bfM (  u )= f &  \quad \text{in $\Omega$}
\end{align}
where $\Div$ and $\text{div}$ denote the divergence of a tensor and a vector 
field, respectively.
We shall first consider a smooth boundary $\Gamma$ with simply supported boundary conditions
\begin{equation}\label{eq:plate-ss}
u = 0 \quad \text{on $\Gamma$}, \qquad M_{nn}(u) = 0 \quad \text{on $\Gamma$}
\end{equation}
where $M_{ab} ={\boldsymbol a}\cdot {\boldsymbol M}\cdot {\boldsymbol b}$ for ${\boldsymbol a},{\boldsymbol b} \in \bbR^2$.
Defining the tangent vector on the boundary as $\bft = (n_2,-n_1)$, multiplying by a test function $v$ and using repeated integration by parts we find that
\begin{align}
(\text{div}\, \Div\, \bfM(  u ),v) ={}& (\bfM( u),\bfkappa(v)) -\langle M_{nn}(u),\partial_n v\rangle_\Gamma
\\
 {}& -\langle M_{nt}(u),\partial_t v\rangle_\Gamma +\langle\bfn\cdot \Div\, \bfM(u) ,v\rangle_\Gamma
 \end{align}
In the case of a smooth boundary we note that
\begin{equation}\label{eq:twist}
\langle M_{nt}(u),\partial_t v\rangle_\Gamma = -\langle \partial_tM_{nt}(u),v\rangle_\Gamma
\end{equation}
and by introducing the Kirchhoff shear force $T :=\bfn\cdot \Div\, \bfM +\partial_tM_{nt}$ we have
\begin{equation}
(\text{div}\, \Div\, \bfM (  u ),v) = (\bfM( u),\bfkappa(v)) -\langle M_{nn}(u),\partial_n v\rangle_\Gamma
+\langle T(u), v\rangle_\Gamma 
 \end{equation}
Taking into account the boundary conditions, the variational problem thus takes the form: 
find 
\[ u \in V = \{ v \in H^2(\Omega) : \text{$v=0$ on $\Gamma$}\}\] such that
\begin{equation}\label{eq:varform}
(\bfM(u), \bfkappa(v) )= (f,v) \quad \forall v \in V
\end{equation}

We will next consider the Signorini condition $u\geq g$ on $\Gamma$, which corresponds to a case where the plate boundary rests on a rigid foundation but is not fixed to it. Introducing a multiplier representing $T(u)$ we have that
\begin{align}
\text{div}\, \Div\,\bfM( u )={}& f   \quad \text{in $\Omega$}\\
M_{nn}(u) ={}& 0 \quad \text{on $\Gamma$}\\
T(u)+\lambda ={}& 0 \quad \text{on $\Gamma$}\\
u-g \geq {}& 0\quad \text{on $\Gamma$}\\
\lambda \leq {}& 0\quad \text{on $\Gamma$}\\
\lambda(u-g) = {}& 0\quad \text{on $\Gamma$}
\end{align}
In this case, the Kuhn--Tucker conditions can be rewritten
\begin{equation}
\lambda = -{\gamma}[g-u -\gamma^{-1} \lambda]_+
\end{equation}
\begin{rem}[Handling polygonal domains]
In the case of a domain with piecewise smooth boundaries, so called Kirchhoff corner forces occur in corner points \cite[Chapter 5.5]{Sl06}.
This case was considered by Nazarov et al. \cite{NaStSw12} but with an alternative formulation (the biharmonic operator, leading to quite different boundary conditions). We here assume that $\Gamma$ consists of smooth connected parts $\Gamma_i$ with corner intersections at $\bfx_i$. Now (\ref{eq:twist}) has to be modified as follows:
\begin{equation}\label{eq:twist2}
\langle M_{nt}(u),\partial_t v\rangle_\Gamma = -\langle \partial_tM_{nt}(u),v\rangle_\Gamma +\sum_i \left(M_{nt}^-(u(\bfx_i)-M_{nt}^+(u(\bfx_i))\right)v(\bfx_i)
\end{equation}
where $M_{nt}^{\pm}(u(\bfx_i)) = \lim_{\epsilon\downarrow 0}M_{nt}(u(x_i\pm\epsilon,y_i\pm\epsilon))$,
giving rise to (virtual work of) point forces in the corners. Unlike the Kirchhoff shear forces, the corner forces are present whether there is contact or not, and are implemented as contributions to the stiffness matrix.
\end{rem}

\section{Finite element method}

We will use $C^1$--continuous element on meshes $\mcK_h$ made up of rectangles. On each element $K\in \mcK_h$ we let $Q_3$ denote the outer product of cubic polynomials:
\[
Q_3 =\left\{ p(x,y): \, p(x,y) =\sum_{0\leq i,j \leq 3} c_{ij} x^i y^j\right\}
\]
where $c_{ij}$ are constants.
The approximation space associated with the Bogner-Fox-Schmit (BFS) element first proposed in \cite{BoFoSc65} is defined by
\begin{equation}
V_{h} = \left\{v\in C^1(\Omega): v\vert_K\in Q_3, \; \forall K\in\mcK_h\right\}
\end{equation}
The shape functions on the BFS element are then made up of outer products of cubic splines, typically used for beam problems. We refer to Zhang \cite{Zh10}
for further details on this approximation. Though this element might seem limited in view of it only being defined on rectangular meshes, the recent
CutFEM for BFS \cite{burman2019cut} extends its use to arbitrary geometries.

In analogy with (\ref{eq:contact_Nitsche}) we now pose the following discrete problem: find $(u_h,\lambda_h)\in V_h\times \Lambda_h$, $\Lambda_h$ to be chosen, such that
\begin{equation}
(\bfM(u_h),\bfkappa(v)) -\lga{{\gamma}\left[g-{u_h}-\gamma^{-1} \lambda_h\right]_+}{v} = ({f},{v}) \quad \forall {v}\in V_h
\end{equation}
and
\begin{equation}\label{aug2h}
-\lga{{\gamma}\left[g-{u_h}-\gamma^{-1}\lambda_h\right]_+}{\gamma^{-1}\mu} -\lga{\gamma^{-1} \lambda_h}{ \mu} =0   \quad \forall \mu\in \Lambda_h 
\end{equation}
We next consider replacing $\lambda_h$ following the ideas of Section \ref{sec:1}. To this end, we formally set $\lambda_h = -T(u_h)$ and $\mu = - T(v)$ to obtain the problem of finding $u_h\in V_h$ such that
\begin{equation}
(\bfM(u_h),\bfkappa(v)) -\lga{{\gamma}\left[g-\psi(u_h) \right]_+}{\psi(v)}-\lga{\gamma^{-1} T(u_h)}{T(v)}= ({f},{v})
\end{equation}
for all ${v}\in V_h$, where $\psi(w):={w}-\gamma^{-1} T(w)$. Setting
now $\gamma =\gamma_0/h^3$ stability, existence and uniqueness of discrete
solution can be shown combining the results of \cite{HaLa02}
and \cite{BuHaLa19}. We leave the details to a forthcoming publication.

\section{Numerical results}

We consider a quadratic plate $(0,1)\times(0,1)$ of thickness $t=0.1$ and with moduli of elasticity $E=100$, $\nu=0.5$. This plate is loaded by a point force
of unit strength. The free parameter was chosen as $\gamma_1 = 10^4D$. The maximum displacement on the boundary is set to $g=0$.

\subsection{Point load in the center of the plate}

\begin{figure}[hb]
%\sidecaption
% For example, with the graphicx style use
\includegraphics[scale=.3]{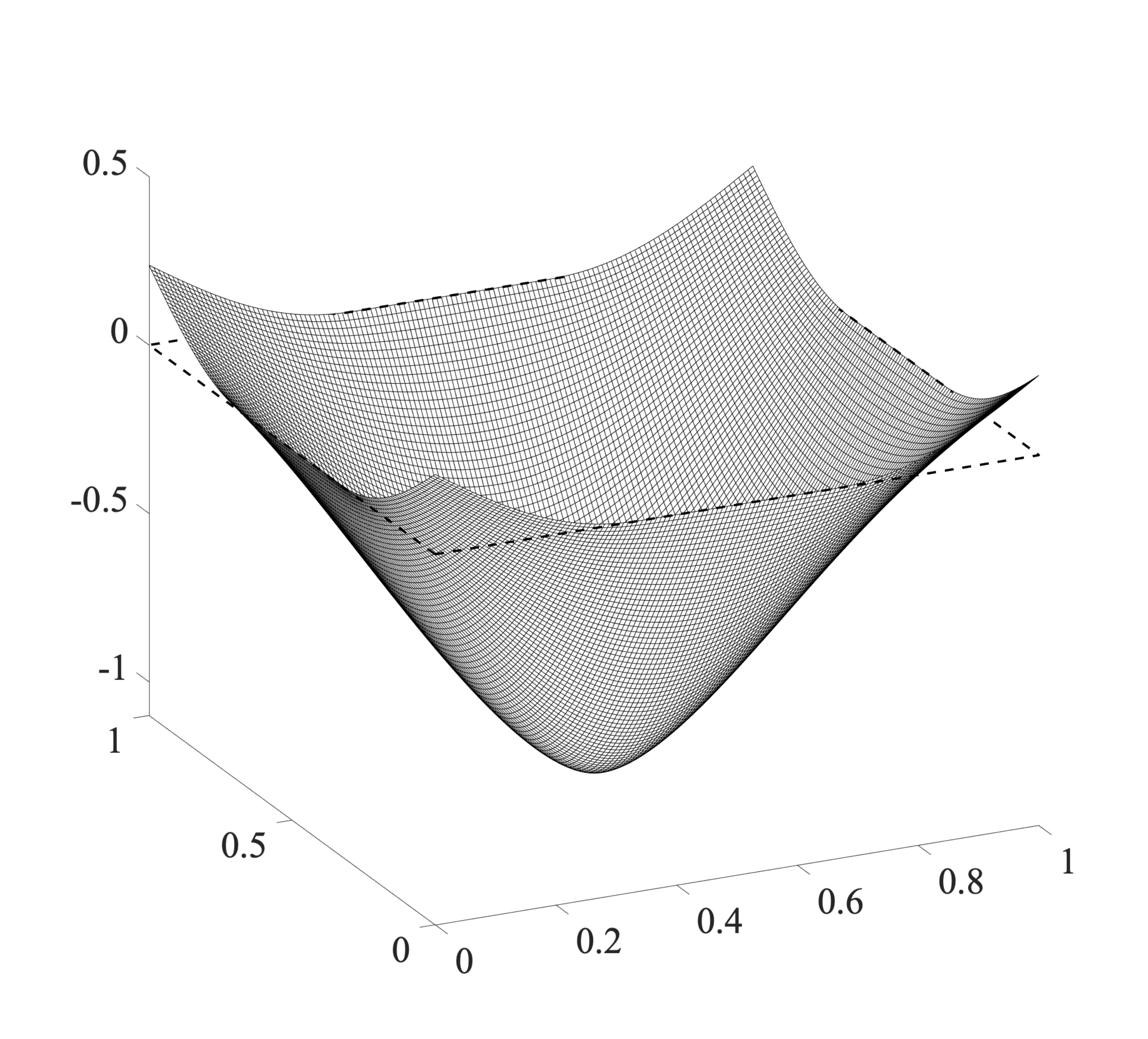}
%
% If no graphics program available, insert a blank space i.e. use
%\picplace{5cm}{2cm} % Give the correct figure height and width in cm
%
\caption{Elevation of the solution on the finest mesh in a sequence. Point load at $(1/2,1/2)$.}
\label{fig:1}       % Give a unique label
\end{figure}

\begin{figure}[hb]
%\sidecaption
% For example, with the graphicx style use
\includegraphics[scale=.16]{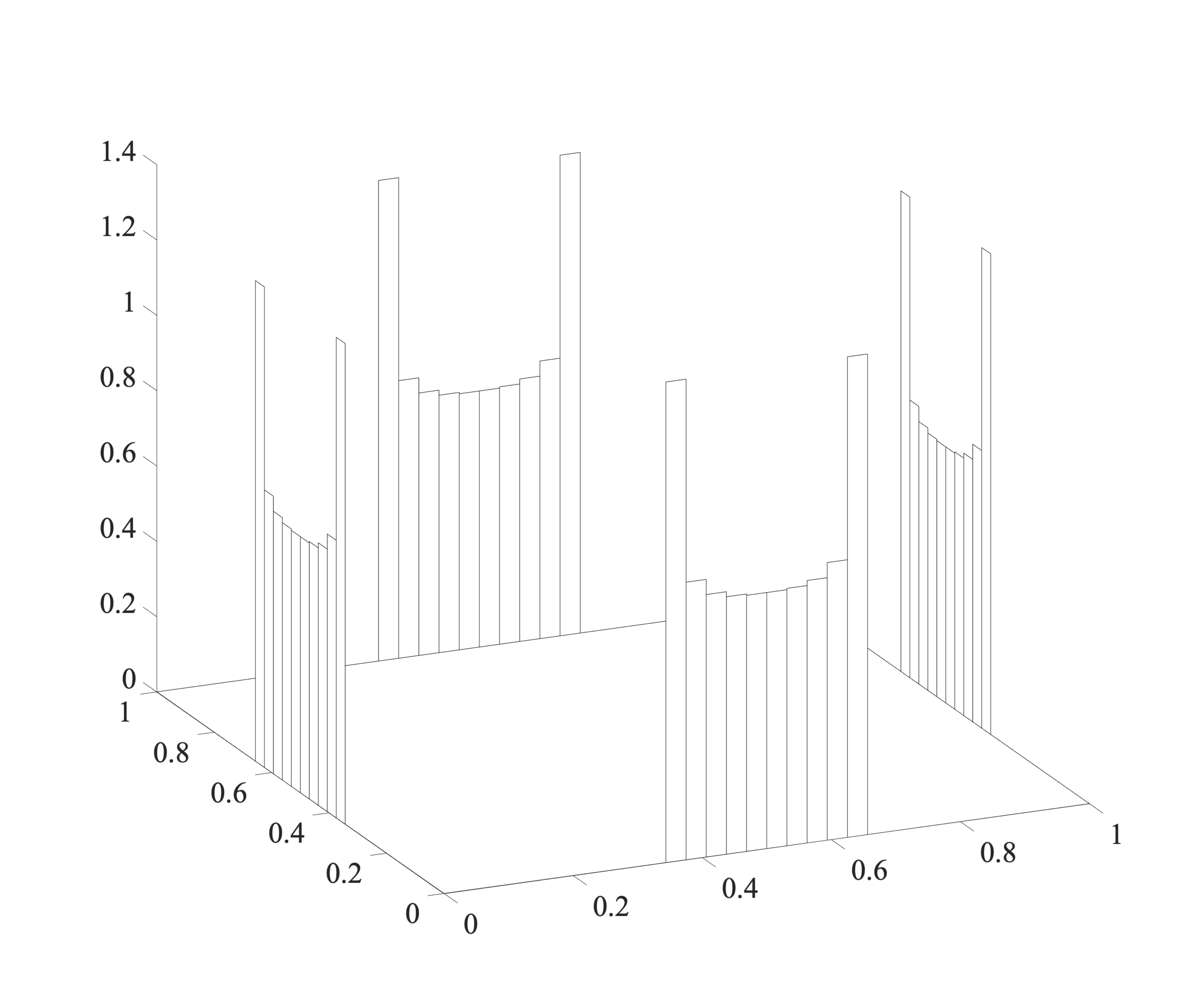}\newline\includegraphics[scale=.22]{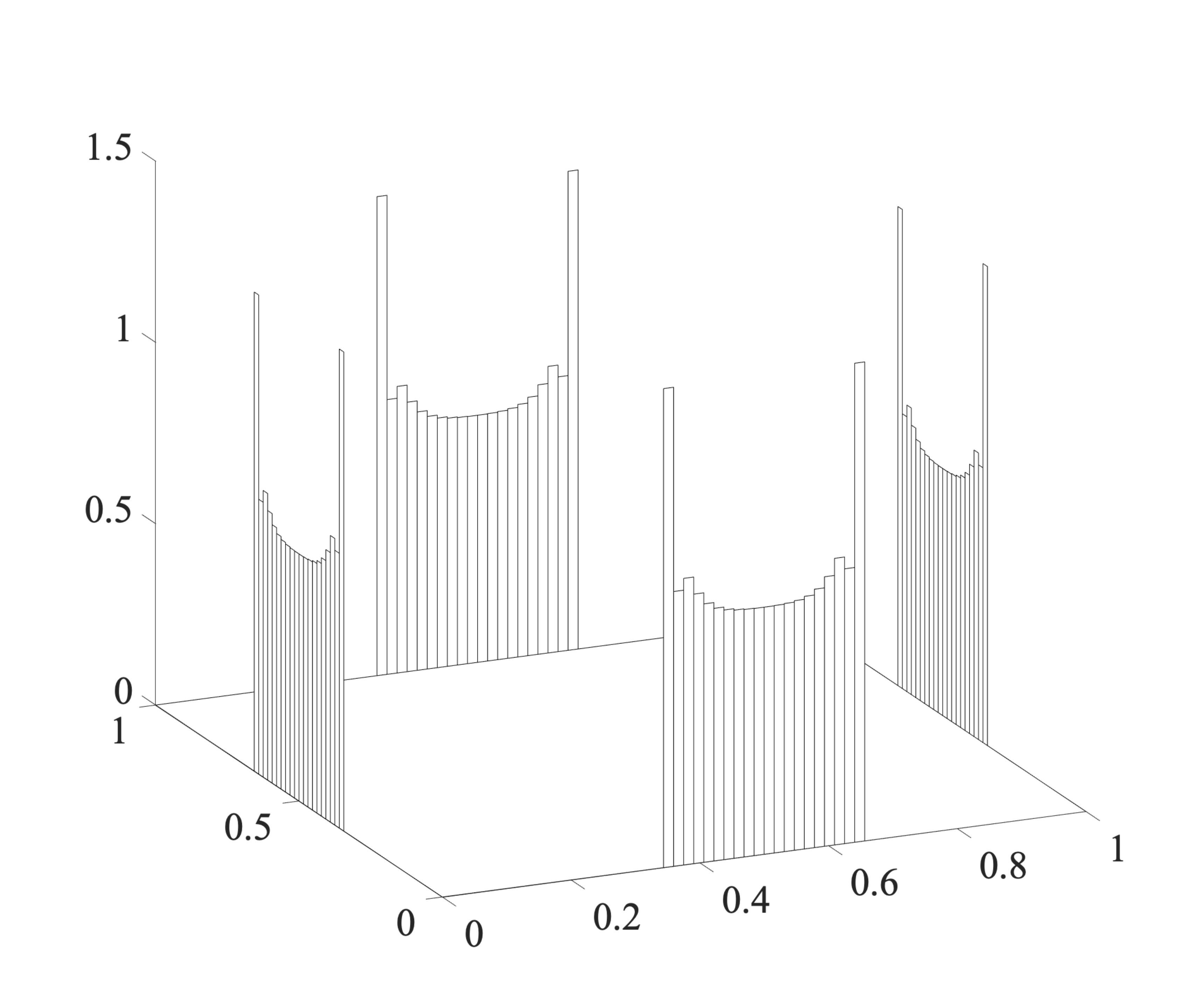}\newline\includegraphics[scale=.22]{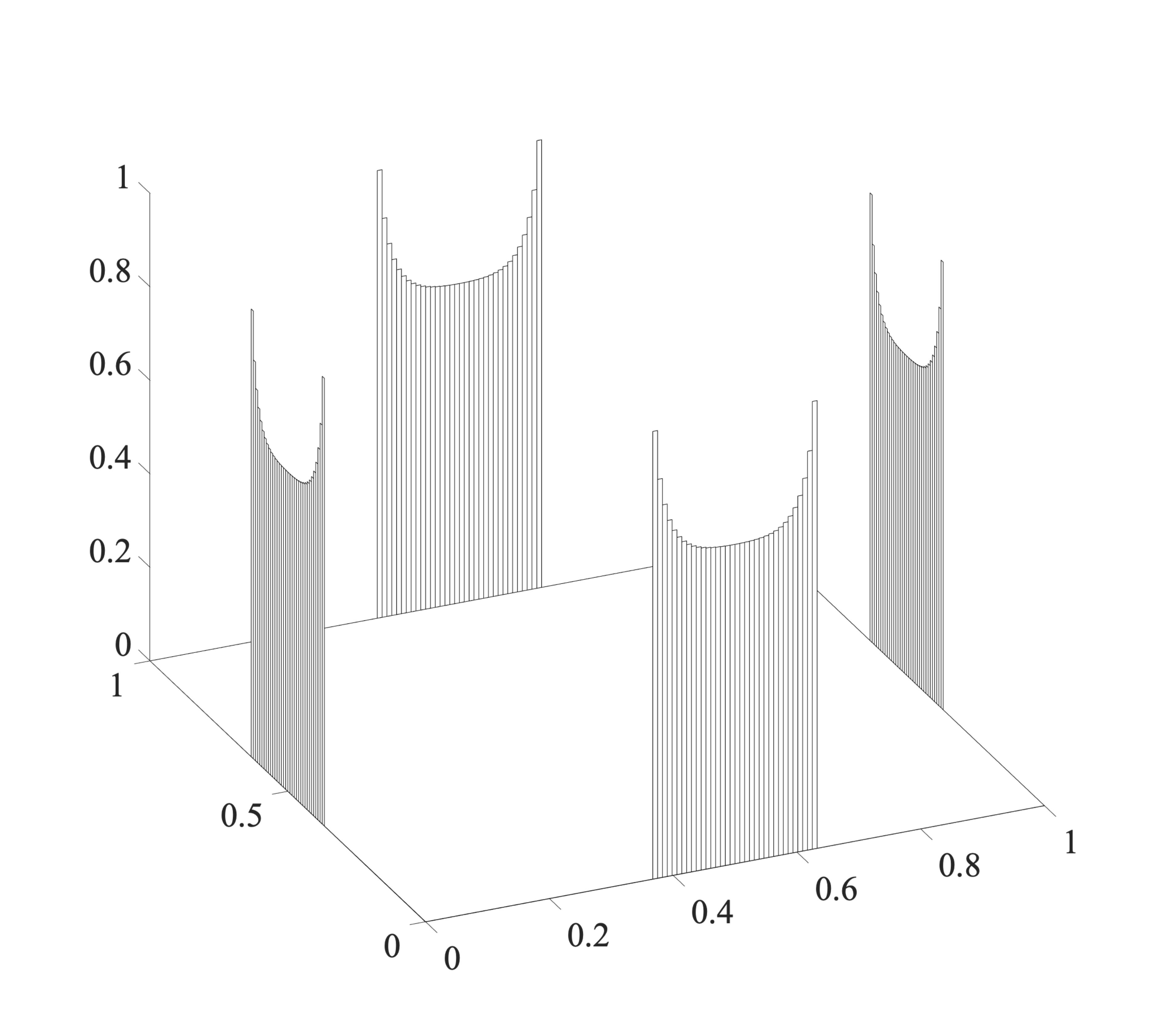}
%
% If no graphics program available, insert a blank space i.e. use
%\picplace{5cm}{2cm} % Give the correct figure height and width in cm
%
\caption{Kirchhoff shear forces in the contact zone on consecutively refined meshes. }
\label{fig:2}       % Give a unique label
\end{figure}

We load the plate with a unit point load at the center. In Fig. \ref{fig:1} we show the computed displacement field with the Signorini boundary indicated by a dotted line. In Fig \ref{fig:2} we show the computed Kirchhoff shear force in the contact zone (evaluated at the midpoint of each element side) on a sequence of uniformly refined meshes. We note the symmetry of the solution.

\subsection{Point load at $(3/4,3/4)$}
\begin{figure}[hb]
%\sidecaption
% For example, with the graphicx style use
\includegraphics[scale=.3]{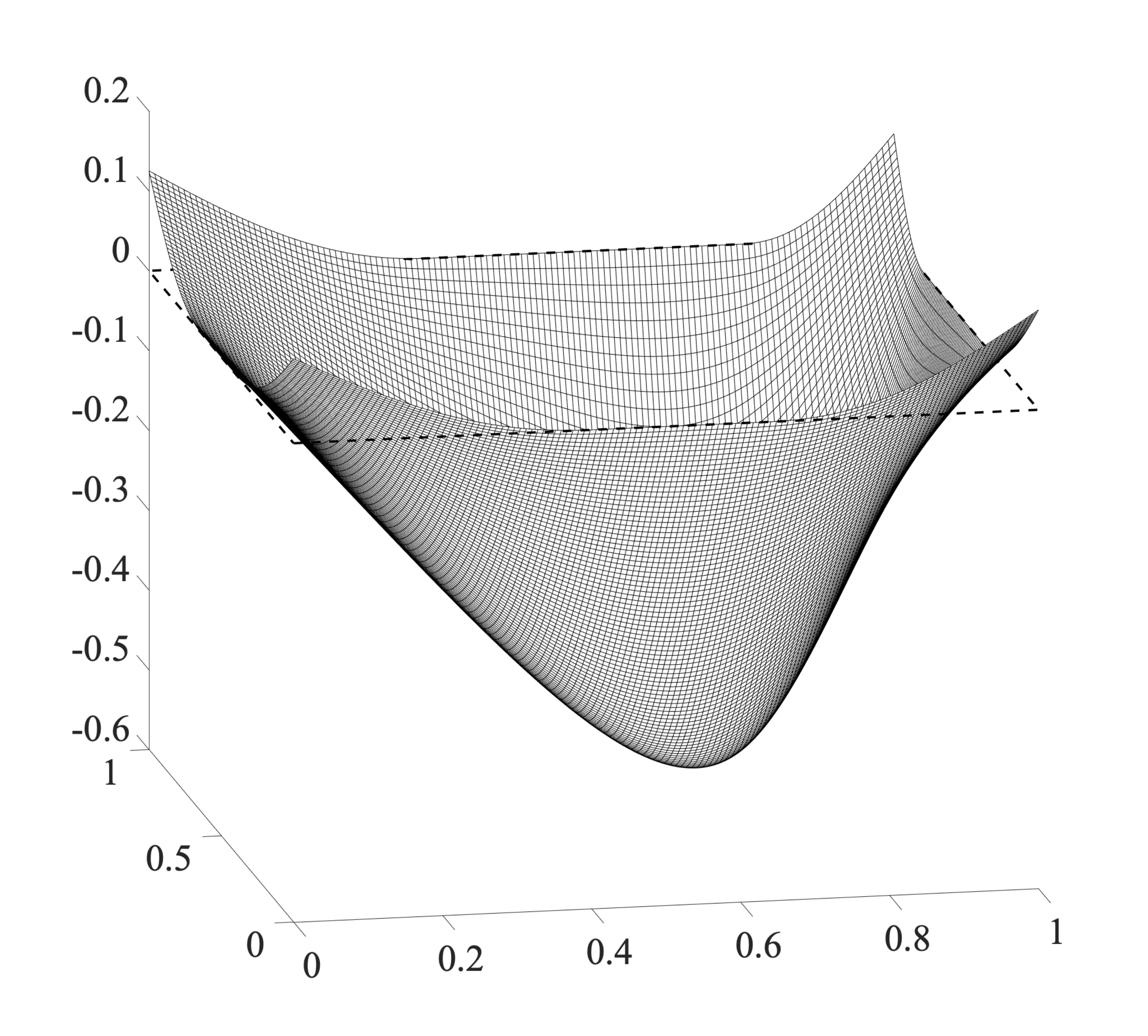}
%
% If no graphics program available, insert a blank space i.e. use
%\picplace{5cm}{2cm} % Give the correct figure height and width in cm
%
\caption{Elevation of the solution on the finest mesh in a sequence. Point load at $(3/4,3/4)$.}
\label{fig:3}       % Give a unique label
\end{figure}

\begin{figure}[hb]
%\sidecaption
% For example, with the graphicx style use
\includegraphics[scale=.17]{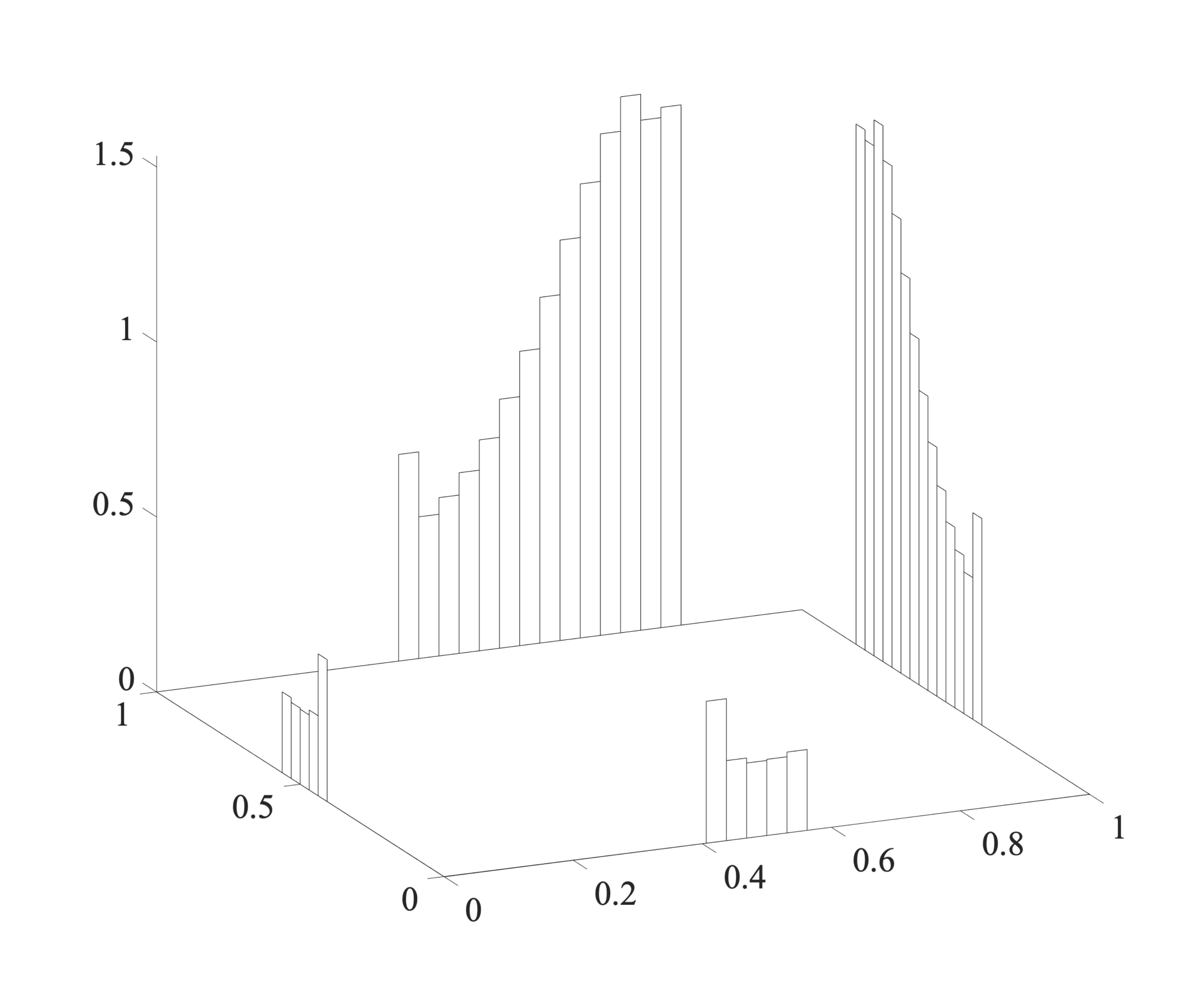}\newline\includegraphics[scale=.18]{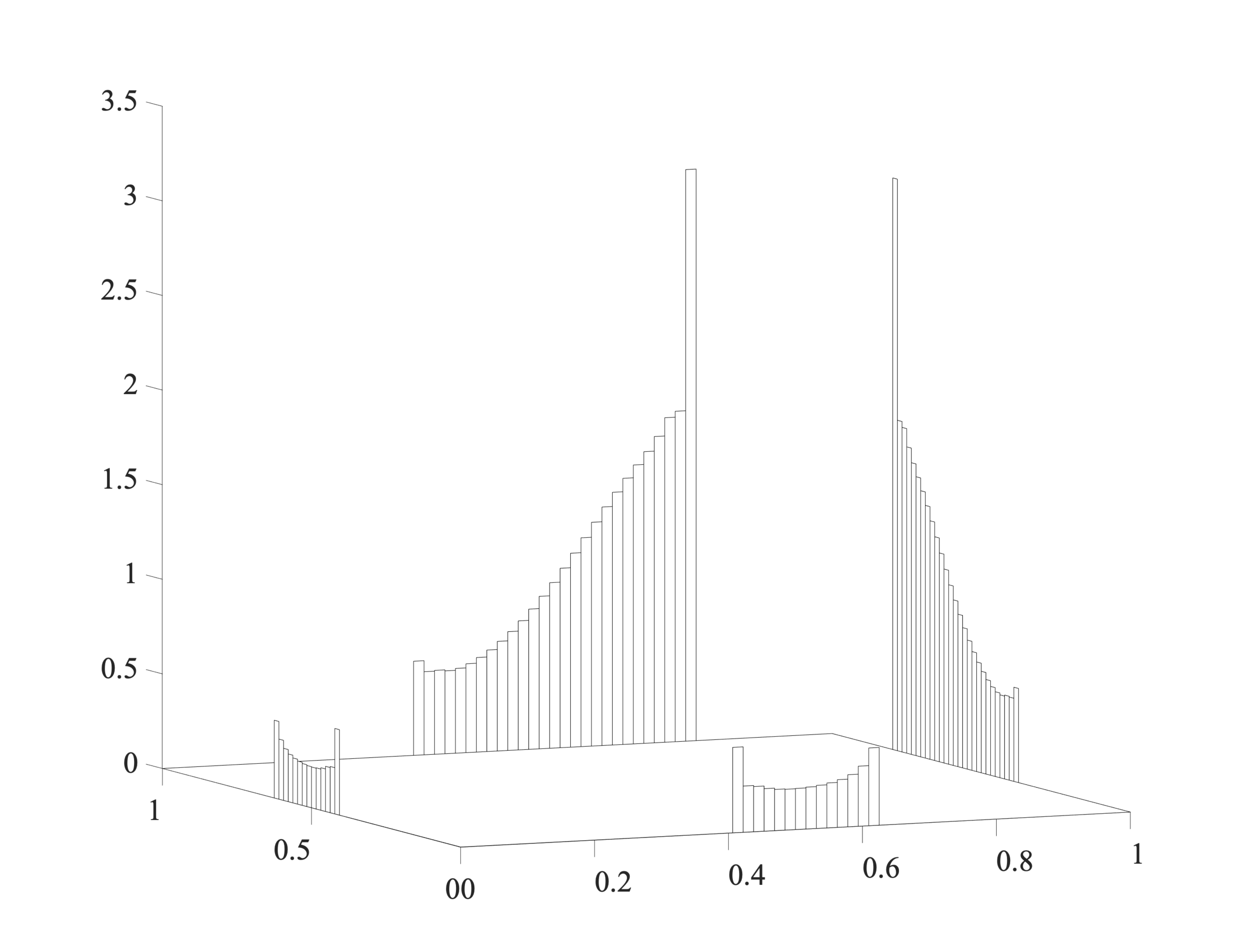}\newline\includegraphics[scale=.22]{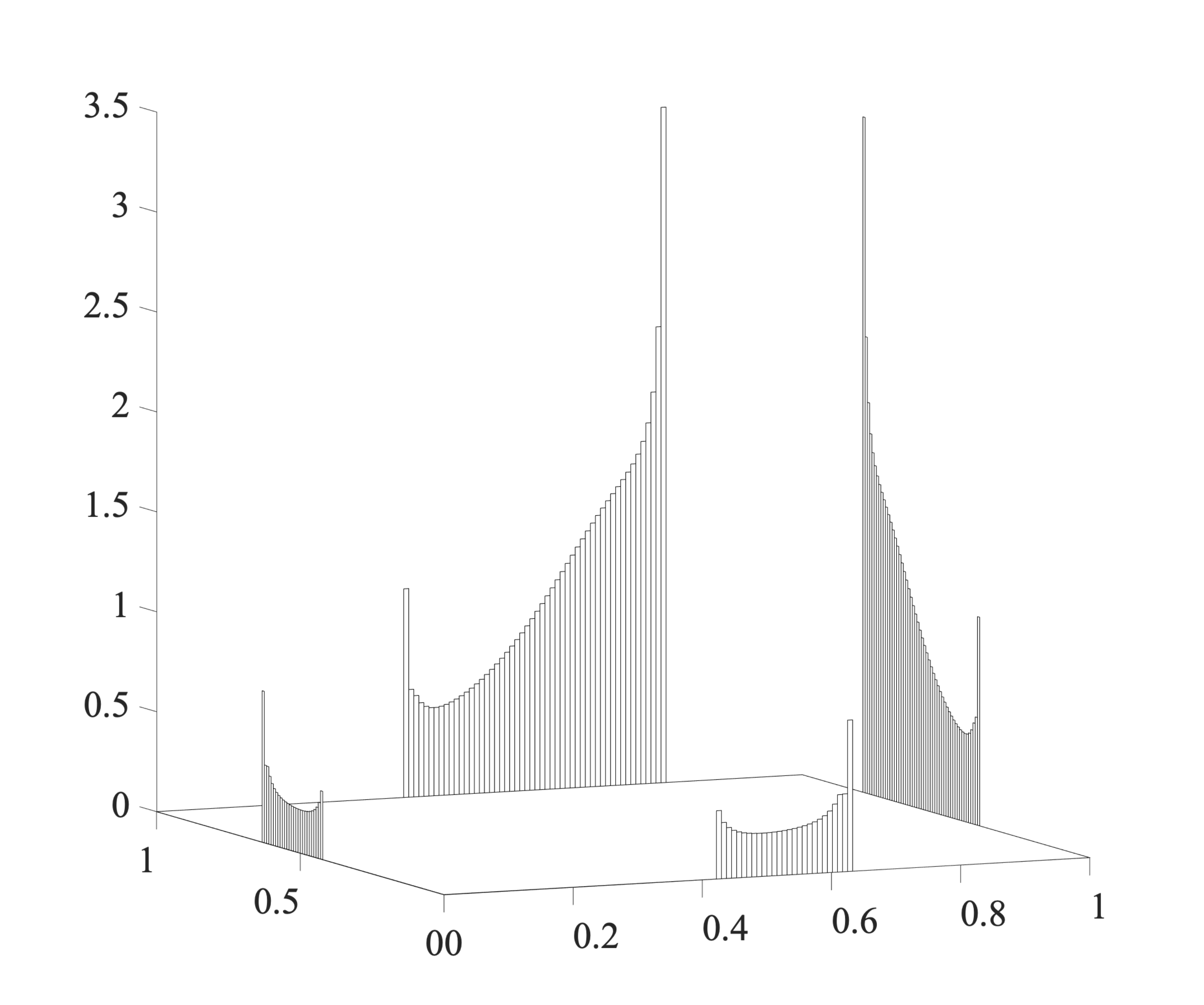}
%
% If no graphics program available, insert a blank space i.e. use
%\picplace{5cm}{2cm} % Give the correct figure height and width in cm
%
\caption{Kirchhoff shear forces in the contact zone on consecutively refined meshes. }
\label{fig:4}       % Give a unique label
\end{figure}

The same plate is now loaded with unit point load at $(3/4,3/4)$. In Fig. \ref{fig:3} we show the computed displacement field, again with the Signorini boundary indicated by a dotted line. In Fig \ref{fig:4} we show the corresponding Kirchhoff shear force in the contact zone. We note the elevation of the shear force close to the first point of contact, similar, but more pronounced, to Fig. \ref{fig:2}. 

\section*{Acknowledgments}
This research was supported in part by the Swedish Foundation
for Strategic Research Grant No.\ AM13-0029, the Swedish Research
Council Grants Nos.\  2017-03911 and 2018-05262 , and the Swedish
Research Programme Essence. Erik Burman was partially supported by 
the grant: EP/P01576X/1.


\end{document}